\documentclass[a4paper,12pt]{article}

\usepackage[latin1]{inputenc}
\usepackage[UKenglish]{babel}
\usepackage[UKenglish]{isodate}

\usepackage{graphicx,amssymb,amsmath,enumerate}
\usepackage{amsfonts,titling,setspace}
\usepackage{amsthm,titlesec,units,makeidx}

\usepackage{tikz}
\usetikzlibrary{arrows.meta}

\newtheorem{theorem}{Theorem}
\newtheorem{corollary}[theorem]{Corollary}

\title{A sharp upper bound for the harmonious total chromatic number of  graphs and multigraphs}

\author{M. Abreu \protect\footnotemark[1]
\and J. B. Gauci \protect\footnotemark[2]
\and D. Mattiolo \protect\footnotemark[3]
\and G. Mazzuoccolo \protect\footnotemark[4]
\and F. Romaniello \protect\footnotemark[5]
\and C. Rubio-Montiel \protect\footnotemark[6]
\and T. Traetta \protect\footnotemark[7]}

\date{\today}

\begin{document}
\maketitle

\def\thefootnote{\fnsymbol{footnote}}
\footnotetext[1]{Universit\`a degli Studi della Basilicata, Via dell'Ateneo Lucano 10, 85100 Potenza, Italy.}
\footnotetext[2]{Department of Mathematics, Faculty of Science, University of Malta, Msida, Malta.}
\footnotetext[3]{Department of Computer Science, KU Leuven Kulak, 8500 Kortrijk, Belgium.}
\footnotetext[4]{Dipartimento di Scienze Fisiche, Informatiche e Matematiche, Universit\`a di Modena e Reggio Emilia, Via Campi 213/b 41125 Modena, Italy.}
\footnotetext[5]{Dipartimento di Matematica ``Giuseppe Peano",  Universit\`a di Torino, via Carlo Alberto 10, 10123 Torino, Italy.}
\footnotetext[6]{Divisi\'on de Matem\'aticas e Ingenier\'ia, FES Acatl\'an, Universidad Nacional Aut\'onoma de M\'exico, Naucalpan, Mexico.}
\footnotetext[7]{DICATAM - Sez. Matematica, Universit\`a degli Studi di Brescia, Via Valotti 9, I-25123 Brescia, Italy.}

\footnotetext[1]{Supported by the group GNSAGA of INdAM.}
\footnotetext[3]{Supported by a Postdoctoral Fellowship of the Research Foundation Flanders (FWO), project number 1268323N.}
\renewcommand{\thefootnote}{\arabic{footnote}}

\begin{abstract}
A proper total colouring of a graph $G$ is called harmonious if it has the further property that when replacing each unordered pair of incident vertices and edges with their colours, then no pair of colours appears twice. The smallest number of colours for it to exist is called the harmonious total chromatic number of $G$, denoted by $h_t(G)$. Here, we give a general upper bound for $h_t(G)$ in terms of the order $n$ of $G$. Our two main results are obvious consequences of the computation of the harmonious total chromatic number of the complete graph $K_n$ and of the complete multigraph $\lambda K_n$, where $\lambda$ is the number of edges joining each pair of vertices of $K_n$. In particular, Araujo-Pardo et al. have recently shown that  $\frac{3}{2}n\leq h_t(K_n) \leq \frac{5}{3}n +\theta(1)$. In this paper, we prove that $h_t(K_{n})=\left\lceil \frac{3}{2}n \right\rceil$ except for $h_t(K_{1})=1$ and $h_t(K_{4})=7$; therefore, $h_t(G) \le \left\lceil \frac{3}{2}n \right\rceil$, for every graph $G$ on $n>4$ vertices.
Finally, we extend such a result to the harmonious total chromatic number of the complete multigraph $\lambda K_n$ and as a consequence show that $h_t(\mathcal{G})\leq (\lambda-1)(2\left\lceil\frac{n}{2}\right\rceil-1)+\left\lceil\frac{3n}{2}\right\rceil$ for $n>4$, where $\mathcal{G}$ is a multigraph such that  $\lambda$ is the maximum number of edges between any two vertices.
\end{abstract}
\textbf{Math. Subj. Class.}: 05C15, 05C70.

\section{Introduction}

A harmonious colouring of a graph $G$ is a proper vertex colouring such that every pair of colours appears on the end vertices of at most one edge of $G$. The harmonious chromatic number, denoted by $h(G)$, is the smallest number of colours required for a harmonious colouring of $G$. This parameter was introduced by Miller and Pritikin in \cite{MP_harmonious_cn}. It is a small variation of the parameter defined independently by Hopcroft and Krishnamoorthy in \cite{HK_Harm_coloring_of_graphs} and by Frank, Harary, and Plantholt in \cite{FHP_l_distinguishing_number} where vertex colourings are not required to be proper.

Note that a vertex colouring where all pairs of vertices have different colours is harmonious. Thus, for every graph $G$, $h(G)\le|V(G)|$. If $G$ has diameter $2$, then equality holds. Indeed, note that any two vertices with a common neighbour must be coloured differently in a harmonious colouring. On the other hand, the problem of finding $h(G)$ for graphs with diameter at least three is quite difficult. It is proved in \cite{HK_Harm_coloring_of_graphs} that determining $h(G)$ is NP-complete (in particular, in the appendix of \cite{HK_Harm_coloring_of_graphs} a short and elegant proof by David S.\ Johnson is presented). Results on $h(G)$ can be nevertheless achieved by considering families of graphs with specific structural properties, see for example \cite{AAEEJR, AKK, BMNRM, EMcD, Mitchem}, and interesting upper bounds for general graphs can be found for example in \cite{Lu, VVK}.

In \cite{AMOR}, the authors give a lower bound and an upper bound for the harmonious chromatic number of the Levi graph of the complete graph $K_n$ on $n$ vertices, stated here in Theorem \ref{thm:knownbounds} hereunder).
The \emph{Levi graph} of a graph $G$ is the incidence graph of $G$. Recall that a total colouring of a graph is a colouring of both vertices and edges, where pairs of adjacent vertices, pairs of adjacent edges, and pairs of incident vertices and edges receive different colours. A harmonious colouring of the Levi graph of $G$ is then a proper total colouring of $G$ with the additional property that, if we write the colours occurring on a vertex and each edge incident to it as unordered pairs of colours, then each pair would occur at most once. We refer to such a total colouring as a \emph{harmonious total colouring}. Hence, we can refer to the harmonious chromatic number of the Levi graph of a graph $G$ as the harmonious total chromatic number of $G$, denoted by $h_t(G)$.

In this paper we first furnish a general upper bound for $h_t(G)$ in terms of the order $n$ of the graph $G$, and we prove it to be optimal in all cases. 
It is straightforward that $h_t(H) \leq h_t(G)$ for every subgraph $H$ of a graph $G$. Hence, the computation of $h_t(K_n)$ proves a general upper bound for the harmonious total chromatic number of a graph of order $n$. 

In this direction, the following result was proved in \cite{AMOR}.

\begin{theorem}\cite{AMOR}\label{thm:knownbounds}
Let $K_n$ be the complete graph on $n\geq 2$ vertices. Then, the following holds
$$\frac{3}{2}n\leq h_t(K_n) \leq \frac{5}{3}n + \theta(1).$$
\end{theorem}

Here we improve such a result by giving the exact value for $h_t(K_n)$, and our result is the following.

\begin{theorem}\label{theo_main} 
The harmonious total chromatic number $h_t(K_n)$ of the complete graph $K_n$ equals $\left\lceil\frac{3}{2}n\right\rceil$ for $n\not = 1,4$, while $h_t(K_1)=1$ and $h_t(K_4)=7$.
\end{theorem}

\begin{figure}[!ht]
\begin{center}
\begin{tabular}{c}
\begin{tikzpicture}[scale=0.7]
\begin{scope}[every node/.style={circle,thick,draw}]
    \node (1) at (2,5) {$\infty_1$};
    \node (2) at (5,2) {$0_1$};
    \node (3) at (5,-2) {$1_1$};
    \node (4) at (2,-5) {$2_1$};
    \node (5) at (-2,-5) {$2_0$};
    \node (6) at (-5,-2) {$1_0$};
    \node (7) at (-5,2) {$0_0$};
    \node (8) at (-2,5) {$\infty_0$};
\end{scope}

\begin{scope}[
              every node/.style={fill=white,circle},
              every edge/.style={draw=black,very thick}]
    \path [-] (2) edge node[pos=0.56] {$0$} (8);
    \path [-] (3) edge node[pos=0.47] {$2$} (5);
    \path [-] (4) edge node[pos=0.44] {$0$} (6);
    \path [-] (5) edge node[pos=0.47] {$2_1$} (7);
    \path [-] (6) edge node[pos=0.44] {$0_0$} (8);
    \path [-] (1) edge node {$2_1$} (2);
    \path [-] (1) edge node[pos=0.47] {$0_1$} (3);
    \path [-] (1) edge node {$1_1$} (4);
    \path [-] (1) edge node[pos=0.4] {$0$} (5);
    \path [-] (1) edge node {$1$} (6);
    \path [-] (1) edge node[pos=0.53] {$2$} (7);
    \path [-] (1) edge node {$\infty$} (8);
    \path [-] (2) edge node {$0_0$} (3);
    \path [-] (2) edge node[pos=0.44] {$1_0$} (4);
    \path [-] (2) edge node {$1$} (5);
    \path [-] (2) edge node[pos=0.4] {$2$} (6);
    \path [-] (2) edge node {$\infty$} (7);
    \path [-] (3) edge node {$\infty_0$} (4);
    \path [-] (3) edge node {$\infty$} (6);
    \path [-] (3) edge node[pos=0.4] {$0$} (7);
    \path [-] (3) edge node {$1$} (8);
    \path [-] (4) edge node {$\infty$} (5);
    \path [-] (4) edge node {$1$} (7);
    \path [-] (4) edge node[pos=0.4] {$2$} (8);
    \path [-] (5) edge node {$1_1$} (6);
    \path [-] (5) edge node {$1_0$} (8);
    \path [-] (6) edge node {$\infty_1$} (7);
    \path [-] (7) edge node {$2_0$} (8);
\end{scope}
\end{tikzpicture}\\
  $K_8$\\
\end{tabular}
\caption{Harmonious total colourings of $K_8$ with 12 colours.}\label{Fig1}
\end{center}
\end{figure}
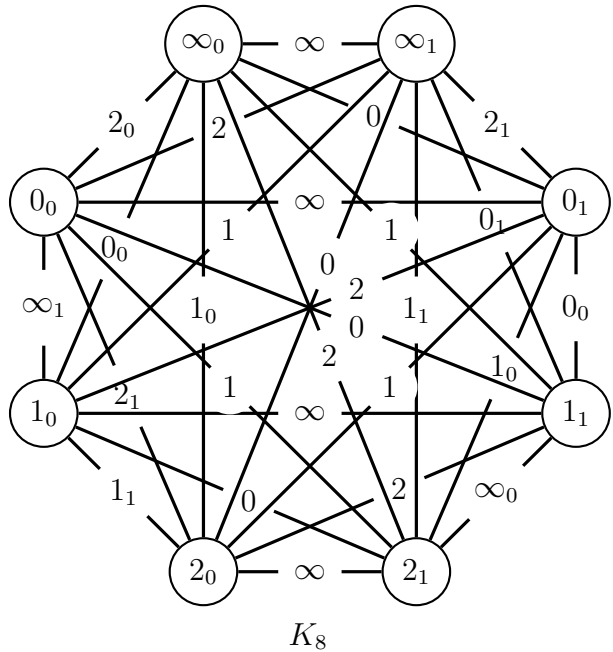

An illustration of a harmonious total colouring of $K_8$ obtained as described in Section 2.2 is shown in Figure \ref{Fig1}. 

As already remarked, our first main result is  a natural consequence of Theorem \ref{theo_main}.

\begin{corollary}
Let $G$ be a graph of order $n$ different from $K_1$ and $K_4$. Then,
$$ h_t(G) \leq 
\left\lceil\frac{3}{2}n\right\rceil.$$
\end{corollary}

The next section is devoted to the proof of Theorem \ref{theo_main}.

Finally, in Section \ref{thirdSec}, we consider the \emph{complete multigraph} $\lambda K_{n}$ defined as the multigraph with $n$ vertices and $\lambda$ edges joining each pair of vertices. We prove the following result.

\begin{theorem}\label{theo_multi} 
The harmonious total chromatic number $h_t(\lambda K_{n})$ of the complete multigraph $\lambda K_{n}$ equals $(\lambda-1)(2\left\lceil \frac{n}{2}\right\rceil -1)+\left\lceil \frac{3n}{2}\right\rceil$ for $2\leq n\not= 4$ and $h_t(\lambda K_{4})=3\lambda+4$.
\end{theorem}

Our second main  result is thus a direct consequence of Theorem \ref{theo_multi}.

\begin{corollary}
     Let $\mathcal{G}$ be a multigraph of order $n\geq 2$, $n\neq 4$, such that $\lambda$ is the maximum number of edges between any two vertices. Then, $$h_t(\mathcal{G})\leq (\lambda-1)(2\left\lceil\tfrac{n}{2}\right\rceil-1)+\left\lceil\tfrac{3n}{2}\right\rceil.$$ 
\end{corollary}

\section{Proof of Theorem \ref{theo_main}}

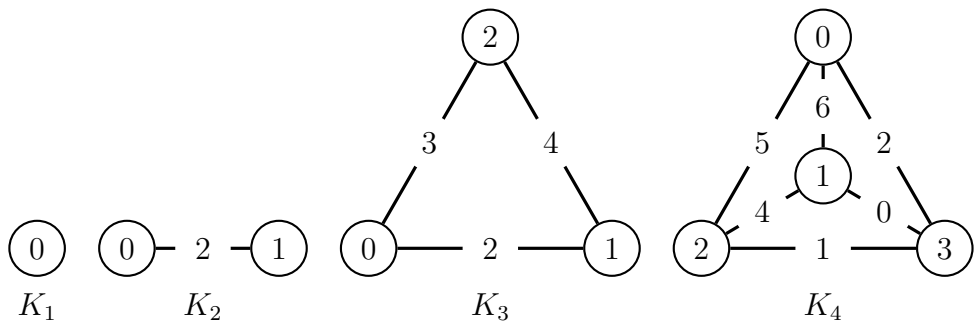
\begin{figure}[!ht]
\begin{center}
\begin{tabular}{cccc}
\begin{tikzpicture}[scale=0.8]

\begin{scope}[every node/.style={circle,thick,draw}]
    \node (0) at (0,0) {0};
\end{scope}

\end{tikzpicture}

&

\begin{tikzpicture}[scale=0.8]

\begin{scope}[every node/.style={circle,thick,draw}]
    \node (0) at (0,0) {0};
    \node (1) at (2.5,0) {1};
\end{scope}

\begin{scope}[
              every node/.style={fill=white,circle},
              every edge/.style={draw=black,very thick}]
    \path [-] (0) edge node {$2$} (1);
\end{scope}
\end{tikzpicture}

&

\begin{tikzpicture}[scale=0.8]

\begin{scope}[every node/.style={circle,thick,draw}]
    \node (0) at (0,0) {0};
    \node (1) at (4,0) {1};
    \node (2) at (2,3.5) {2};
\end{scope}

\begin{scope}[
              every node/.style={fill=white,circle},
              every edge/.style={draw=black,very thick}]
    \path [-] (0) edge node {$2$} (1);
    \path [-] (0) edge node {$3$} (2);
	\path [-] (1) edge node {$4$} (2);
\end{scope}
\end{tikzpicture}

&

\begin{tikzpicture}[scale=0.8]
\begin{scope}[every node/.style={circle,thick,draw}]
    \node (0) at (2,3.5) {0};
    \node (1) at (2,1.2) {1};
    \node (2) at (0,0) {2};
    \node (3) at (4,0) {3};
\end{scope}

\begin{scope}[
              every node/.style={fill=white,circle},
              every edge/.style={draw=black,very thick}]
    \path [-] (0) edge node {$6$} (1);
    \path [-] (0) edge node {$5$} (2);
    \path [-] (0) edge node {$2$} (3);
    \path [-] (1) edge node {$4$} (2);
    \path [-] (1) edge node {$0$} (3);
    \path [-] (2) edge node {$1$} (3);
\end{scope}
\end{tikzpicture}\\
    $K_{1}$ & $K_2$ & $K_3$ & $K_4$\\
\end{tabular}
\caption{Harmonious total colourings of some small complete graphs.}\label{Fig2}
\end{center}
\end{figure}

Figure \ref{Fig2} displays equality for $n \leq 3$, and, for $n=4$, it  shows that $h_t(K_4)\leq 7$. It is easy to check that 6 colours do not suffice for a harmonious total colouring of $K_4$, so $h_t(K_4)=7$, and, together with the case $n=1$, these are the only exceptions in which $h_t(K_n) > \left\lceil\frac{3}{2}n\right\rceil$.

For the remainder of this proof, we need to distinguish cases according to the congruence of $n>5$ modulo $4$ (or $8$).

\subsection{The case $n\neq 4k$}
In this subsection we consider the case $n \not\equiv 0 {\bmod 4}$ and construct a harmonious total colouring of $K_n$ which uses $\left\lceil\frac{3}{2}n\right\rceil$ colours. More precisely, we prove for all positive values of $k$ that:
\begin{itemize}
    \item $h_t(K_{4k+1})=6k+2$;
    \item $h_t(K_{4k+2})=6k+3$;
    \item $h_t(K_{4k+3})=6k+5$.
\end{itemize}

The general idea behind our constructions is to first colour each of the vertices using a different colour and decompose $K_n$ into factors. The edges of one of the factors are cleverly coloured using the colours of the vertices, and the remaining edges are then coloured using the unused colours.

The proposed construction is divided into two main steps, with the initial step being the same for all three cases. We begin by identifying the vertex-set of $K_{4k+\ell}$, for $\ell\in\{1,2,3\}$, with the elements of the additive group $(\mathbb{Z}_{4k+\ell},+)$ and we denote the edges of the complete graph by the unordered pairs $\{i,j\}$ where $i,j\in \mathbb Z_{4k+\ell}$ such that $i \neq j$. We make use of a set of colours ${\cal C}= \mathbb Z_{4k+\ell} \cup {\cal C}'$, where $\cal C'$ is a set disjoint from $\mathbb Z_{4k+\ell}$. In other words, we are using the elements of the additive group both for the vertices of the graph and for some of the colours.  Set $c(i)=i$ for all $i \in \mathbb Z_{4k+\ell}$, that is, to every vertex we assign a colour equal to its label. For $h\in\{0,1,\ldots,4k\}$, we set $c(\{i+h,j+h\})=h$ (addition modulo $4k+\ell$), for all $i\neq j \in \{1,\ldots,2k\}$ such that $i+j=2k+1$. In other words, we assign colour $0$ to the $k$ edges $\{i,j\}$ with $i\neq j \in \{1,\ldots,2k\}$ such that $i+j=2k+1$, and then colour $h$ to an edge obtained by adding $h$ to both the ends of an edge with colour $0$. Up to now, we have assigned a colour to $(4k+\ell)k$ edges of $K_{4k+\ell}$. More precisely, each colour is used for exactly $k$ edges and we assigned colours to exactly $2k$ edges incident to each vertex.

The second step of our construction aims to colour the remaining edges and we must distinguish between three different cases, as follows:
\begin{itemize}
    \item \textbf{$K_{4k+1}$:} The set of the uncoloured edges is a $2k$-regular spanning subgraph $H_1$ on an odd number of vertices. We complete the colouring $c$ by choosing an arbitrary proper $(2k+1)$-edge-colouring of $H_1$ with $2k+1$ colours in $\cal C'$;
    \item \textbf{$K_{4k+2}$:} The set of the uncoloured edges is a $(2k+1)$-regular spanning subgraph $H_2$ on an even number of vertices. Observe that $H_2$ is the Cayley graph of the group $\mathbb{Z}_{4k+2}$ with generating set $S$ consisting of $\{0,2,\ldots,4k\}$ and the involution $2k+1$. Therefore, by Corollary 2.3.1 of \cite{Stong}, we get that $H_2$ is $(2k+1)$-edge-colourable. Hence, we complete the colouring $c$ by choosing an arbitrary proper $(2k+1)$-edge-colouring of $H_2$ with $2k+1$ colours in $\cal C'$.
     \item \textbf{$K_{4k+3}$:} The set of the uncoloured edges is a $(2k+2)$-regular spanning subgraph $H_3$ on an odd number of vertices. We complete the colouring $c$ by choosing an arbitrary proper $(2k+3)$-edge-colouring of $H_3$ with $2k+3$ colours in $\cal C'$.
\end{itemize}
It is worth remarking that Vizing's Theorem assures the existence of such edge-colourings in the first and third case.

\textbf{Claim:} The total colouring $c$ of $K_{4k+\ell}$ is harmonious.\\
\textit{Proof of the claim.}
We need to prove that the list of all unordered pairs
$\{c(i),c(\{i,j\})\}$ has no repetition. Since no two vertices share the same colour, we only need to prove that there are no $i\neq j$ and $h \neq k$ such that
\begin{equation}\label{eq:1}
c(i)=c(\{h,k\}) \text { and } c(h)=c(\{i,j\}).
\end{equation}
Note that no vertex receives a colour in ${\cal C}'$, then relation (\ref{eq:1}) cannot hold for edges with a colour in ${\cal C}'$.

By the symmetry of $c$, it suffices to prove that the edges incident to the vertex $0$ and the vertices incident to edges having colour $0$ are a partition of $\{1,\dots,4k+\ell\}$. Indeed, this is the case since, by definition of $c$, if $c(\{i,j\})=0$ then $i,j \in \{1,\ldots,2k\}$. On the other hand, the edges $\{i+h,j+h\}$ are coloured $h$, and, without loss of generality, setting $i+h=0$ we get $h=-i$. Now, for each $j\in\{1,\ldots,2k\}$, $i=2k+1-j$, implying that the $2k$ edges $\{0,j-i\}$ get distinct colours $-i\in\{2k+1,\ldots, 4k+\ell\}$. \hfill $\qed$

\subsection{The case $n=4k$}

In this case we subdivide the graph $K_n$ into two equal complete graphs and a bipartite graph that joins them, we use three sets of $\frac{1}{2}n$ colours that will be used both for the vertices and for the edges under an algebraic rule. The original assignment of colours almost works, except for one edge, so we will have to modify the colours of several edges at the end in order to attain the chromatic total colouring with exactly $\frac{3}{2}n$ colours. 

\subsubsection{The case $k$ even}

We let $k=2m$, and show that $K_{8m}$ admits a harmonious total colouring with $12m$ colours.

Set $\mathcal{C}'=\mathbb{Z}_{4m-1}\,\cup\,\{\infty\}$, and for $u\in \mathbb{Z}_2$, let $\mathcal{C}_u=\{x_u\mid x\in \mathcal{C}'\}$  be a copy of $\mathcal{C}'$. Also,
set $V(K_{8m}) = \mathcal{C}_0\, \cup\, \mathcal{C}_1$. We proceed by constructing a proper total $12m$-colouring $c$ of $K_{8m}$ which turns out to be harmonious when $m=1$; otherwise, $c$ is ``almost'' harmonious, in the sense that there are only two pairs
$(x, \{x,y\}), (z, \{z,t\})\in V(K_{8m})\times E(K_{8m})$ with $x\neq z$ such that
$\{c(x), c(\{x,y\})\} = \{c(z), c(\{z,t\})\}$. We then modify $c$ to obtain the desired harmonious colouring for every $m\geq 2$.
We use as a set of colours ${\cal C}= \mathcal{C}_0\,\cup\,\mathcal{C}_1\,\cup\,
\mathcal{C}'$.

Set $c(x_u)=x_u$ for every vertex $x_u$ ($u=0,1$) of $K_{8m}$.  To colour its edges, we start by decomposing $K_{8m}$ into the edge disjoint union of the following graphs:
\begin{enumerate}
  \item $K_{4m,4m}$ is the complete bipartite graph whose parts are $\mathcal{C}_0$ and $\mathcal{C}_1$;
  \item $\Gamma_0$ and $\Gamma_1$  are copies of $K_{4m}$ with vertex-sets $\mathcal{C}_0$ and $\mathcal{C}_1$, respectively.
\end{enumerate}
Consider a $1$-factorization $\mathcal{F}=\{F_y\mid y\in \mathcal{C}'\}$ of $K_{4m,4m}$ and
set $c(e)=y$ for every edge $e\in F_y$ and every $y\in \mathcal{C}'$. It is left to colour the edges of $\Gamma_0$ and $\Gamma_1$.

If $m=1$, set $c(\{\infty_u, (j+1)_u\})= j_u$ for every $j\in\mathbb{Z}_3$ and $u\in\{0,1\}$, and let
\begin{align*}
& c(\{0_0,1_0\})=\infty_1, c(\{1_0,2_0\})=1_1, c(\{2_0,0_0\})=2_1, \\
& c(\{0_1,1_1\})= 0_0, c(\{1_1,2_1\})= \infty_0, c(\{2_1,0_1\})= 1_0.
\end{align*}
One can easily check that $c$ is harmonious (refer to Figure \ref{Fig1}).

For $m\geq 2$, we first colour the edges
incident with
$\infty_u$ and those of the form $\{x_u, y_u\}$ where $x,y\in \mathbb{Z}_{4m-1}$ and $y-x\neq \pm 1$; for $u=0,1$, $i=0, \ldots, m-2$ and $j\in\mathbb{Z}_{4m-1}$, set
\begin{align*}
  & c(\{\infty_u, (2m-1+j)_u\})= j_u,\\
  & c(\{j_u, (4i+j+3)_u\}) = (2i+j+1)_u, \; \text{and}\\
  & c(\{j_u, (4i+j+5)_u\}) = (2i+j+3)_{u+1}.
\end{align*}
The remaining edges, those of the form $\{x_u, (x+1)_u\}$ will be coloured as follows:
\begin{align*}
c(\{j_0, (j+1)_0\})&=
\begin{cases}
  \infty_1 & \text{if $j=0,2,4, \ldots, 4m-4$},\\
  j_1      & \text{if $j=1,3,5, \ldots, 4m-3, 4m-2$};
\end{cases}\\
c(\{j_1, (j+1)_1\})&=
\begin{cases}
  \infty_0 & \text{if $j=1,3,5, \ldots, 4m-3$},\\
  j_0 & \text{if $j=0,2,4, \ldots, 4m-2$}.
\end{cases}
\end{align*}
One can check that $c$ is almost harmonious; more precisely,
\begin{enumerate}
  \item $c$ is a proper total colouring of $K_{8m}$,
  \item $c(\{(4m-2)_0,0_0\}) = (4m-2)_1$,
    and $c(\{(4m-2)_1,0_1\}) = (4m-2)_0$, and
  \item $c$ is harmonious on $K_{8m}$ minus the edge $\{(4m-2)_0,0_0\}$.
\end{enumerate}
By recolouring the edges $\{(4m-2)_0,0_0\}, \{0_0, (4m-3)_0\}, \{(2m-1)_1, (6m-4)_1\}$ and
$\{(4m-3)_0, (8m-6)_0\}$ in the following way, we obtain the desired harmonious colouring of $K_{8m}$:
\[
 c(\{(4m-2)_0,0_0\})= (2m-1)_1, \;\;\;
c(\{0_0, (4m-3)_0\})= (4m-2)_1,
\]
\[
c(\{(2m-1)_1, (6m-4)_1\})= (4m-3)_0\;\;\;
c(\{(4m-3)_0, (8m-6)_0\})= (4m-4)_1.
\]

\subsubsection{The case $k$ odd}

We let $k=2m+1$, and show that $K_{8m+4}$ admits a harmonious total colouring with $12m+6$ colours.

Set $\mathcal{C}'=\mathbb{Z}_{4m+1}\,\cup\,\{\infty\}$, and for $u\in \mathbb{Z}_2$, let $\mathcal{C}_u=\{x_u\mid x\in \mathcal{C}'\}$  be a copy of $\mathcal{C}'$. Also,
set $V(K_{8m+4}) = \mathcal{C}_0\, \cup\, \mathcal{C}_1$. We proceed by constructing a proper total $(12m+6)$-colouring $c$ of $K_{8m+4}$ which turns out to be ``almost'' harmonious, in the sense that there are only two pairs
$(x, \{x,y\}), (z, \{z,t\})\in V(K_{8m+4})\times E(K_{8m+4})$ with $x\neq z$ such that
$\{c(x), c(\{x,y\})\} = \{c(z), c(\{z,t\})\}$. We then modify $c$ to obtain the desired harmonious colouring for every $m\geq 1$.
We use as a set of colours ${\cal C}= \mathcal{C}_0\,\cup\,\mathcal{C}_1\,\cup\,
\mathcal{C}'$.

Set $c(x_u)=x_u$ for every vertex $x_u$ ($u=0,1$) of $K_{8m+4}$.  
To colour its edges, we start by decomposing $K_{8m+4}$ into the edge disjoint union of the following graphs:
\begin{enumerate}
  \item $K_{4m+2,4m+2}$ is the complete bipartite graph whose parts are $\mathcal{C}_0$ and $\mathcal{C}_1$;
  \item $\Gamma_0$ and $\Gamma_1$  are copies of $K_{4m+2}$ with vertex-sets $\mathcal{C}_0$ and $\mathcal{C}_1$, respectively.
\end{enumerate}
Consider a $1$-factorization $\mathcal{F}=\{F_y\mid y\in \mathcal{C}'\}$ of $K_{4m+2,4m+2}$ and
set $c(e)=y$ for every edge $e\in F_y$ and every $y\in \mathcal{C}'$. It is left to colour the edges of $\Gamma_0$ and $\Gamma_1$.

We first colour the edges incident with
$\infty_u$ and those of the form $\{x_u, y_u\}$ where $x,y\in \mathbb{Z}_{4m+1}$ and $y-x\neq \pm 1$; for $u=0,1$, $0\leq i, i' \leq m-1$ with $i' \neq  m-1$, and $j\in\mathbb{Z}_{4m+1}$, set
\begin{align*}
  & c(\{\infty_u, j_u\})= (j+2m)_{u+1},\\
  & c(\{j_u, (4i+j+3)_u\}) = (2i+j+2)_u, \; \text{and}\\
  & c(\{j_u, (4i'+j+5)_u\}) = (2i'+j+3)_{u+1}.
\end{align*}
The remaining edges, those of the form $\{x_u, (x+1)_u\}$ will be coloured as follows:
\begin{align*}
c(\{j_0, (j+1)_0\})&=
\begin{cases}
  \infty_0 & \text{if $j=0,2,4, \ldots, 4m-2$},\\
  j_1      & \text{if $j=1,3,5, \ldots, 4m-1, 4m$};
\end{cases}\\
c(\{j_1, (j+1)_1\})&=
\begin{cases}
  \infty_1 & \text{if $j=1,3,5, \ldots, 4m-1$},\\
  j_0 & \text{if $j=0,2,4, \ldots, 4m$}.
\end{cases}
\end{align*}
One can check that $c$ is almost harmonious; more precisely,
\begin{enumerate}
  \item $c$ is a proper total colouring of $K_{8m+4}$,
  \item $c(\{4m_0,0_0\}) = 4m_1$,
    and $c(\{4m_1,0_1\}) = 4m_0$, and
  \item $c$ is harmonious on $K_{8m+4}$ minus the edge $\{4m_0,0_0\}$ (or $\{4m_1,0_1\}$).
\end{enumerate}
By recolouring, in the following way, the edges 
$\{0_1, 4_1\}, \{\infty_0, 3_0\}, \{4_0, 1_0\}$ when $m=1$, and the edges
$\{4m_0,0_0\}, \{(4m-1)_0, 0_0\}, \{(2m)_0, (2m-2)_0\}$ and $\{(4m-1)_0, (4m-3)_0\}$ when $m\geq 2$, we obtain the desired harmonious total colouring of $K_{8m+4}$: when $m=1$, set $c(\{0_1, 4_1\})=3_0, c(\{\infty_0, 3_0\})=4_0$ and $c(\{4_0, 1_0\})=0_1$, otherwise set\
\[
 c(\{4m_0,0_0\})= (2m)_0, \;\;\;
c(\{0_0, (4m-1)_0\})= (4m)_1,
\]
\[
c(\{(2m)_0, (2m-2)_0\})= (4m-1)_0\;\;\;
c(\{(4m-1)_0, (4m-3)_0\})= (4m-2)_1. \qed
\]

\section{Proof of Theorem \ref{theo_multi}}\label{thirdSec}

The computation of the harmonious total chromatic number $h_t(\lambda K_{n})$ is divided into two parts. The lower bound is a generalization of the lower bound of Theorem \ref{thm:knownbounds} and the upper bound is a consequence of Theorem \ref{theo_main}.

\begin{proof}
To achieve the lower bound, let 
$c : V \cup E \rightarrow \{c_1,\dots, c_{h_t(\mathcal{G})}\}$ be a harmonious total colouring of $\mathcal{G}:=\lambda K_{n}$, where $ V= \{v_1, \dots, v_n\}$ and for each $v_i\in V$, let $c_i = c(v_i)$, 
since vertices are pairwise adjacent.

For each $i\in \{1,\dots,h_t(\mathcal{G})\}$,  let $E_i = \{e\in E : c (e)= c_i\}$. Let $\mathcal{G}^*$ be the partial multigraph of $\mathcal{G}$ such that its vertices are $V$ and its edges are $E^*=\bigcup\limits _{i=1}^{n}E_{i}$. On the one hand, $|E^*|= \sum\limits_{i=1}^{n} |E_i| = \frac{1}{2}\sum\limits_{i=1}^{n} \deg_{\mathcal{G}^*}(v_i)$ where $\deg_{\mathcal{G}^*}(v_i)$ is the number of edges of $\mathcal{G}^*$ incident to $v_i$. On the other hand, for each  $v_i \in V$, if there is a edge of $\mathcal{G}^*$, incident to $v_i$, coloured $c_j$, then $v_j$ is not incident to any edge of $E_i$. Thus, $|E_i| \leq \frac{n-1-\deg_{\mathcal{G}^*}(v_i)}{2}$.

Therefore, $\sum\limits_{i=1}^{n} |E_i|  \leq \sum\limits_{i=1}^{n} \frac{n-1-\deg_{\mathcal{G}^*}(v_i)}{2}$ which implies that $\sum\limits_{i=1}^{n} |E_i| \leq \frac{n(n-1)}{2}-\sum\limits_{i=1}^{n} |E_i|$. Hence, $\sum\limits_{i=1}^{n} |E_i| \leq \lfloor \frac{n(n-1)}{4}\rfloor$ and therefore 
\[\lambda \binom{n}{2} =\sum\limits_{i=1}^{h_t(\mathcal{G})} |E_i|=\sum\limits_{i=1}^{n} |E_i|+\sum\limits_{i=n+1}^{h_t(\mathcal{G})} |E_i|\leq  \left\lfloor \frac{n(n-1)}{4}\right\rfloor+\sum\limits_{i=n+1}^{h_t(\mathcal{G})}|E_i|.\]
Since, for each $i\in \{1,\dots,h_t(\mathcal{G})\}$, $|E_i|\leq \lfloor\frac{n}{2}\rfloor$ it follows that

\[\lambda\binom{n}{2}\leq\left\lfloor \frac{n(n-1)}{4}\right\rfloor +(h_{t}(\mathcal{G})-n)\left\lfloor \frac{n}{2}\right\rfloor,\]

\[\left\lceil \frac{\lambda\binom{n}{2}-\left\lfloor \frac{n(n-1)}{4}\right\rfloor }{\left\lfloor \frac{n}{2}\right\rfloor }\right\rceil +n\leq h_{t}(\mathcal{G}).\]

After some simple calculations, we get

\[h_{t}(\mathcal{G})\geq\begin{cases}
\begin{array}{c}
(\lambda-1)(n-1)+\frac{3n}{2}\\
(\lambda-1)n+\frac{3n+1}{2}
\end{array} & \begin{array}{c}
\textrm{if }n\textrm{ is even},\\
\textrm{if }n\textrm{ is odd}.
\end{array}\end{cases}\]

Thus we can conclude that $(\lambda-1)(2\left\lceil \frac{n}{2}\right\rceil -1)+\left\lceil \frac{3n}{2}\right\rceil\leq h_t(\lambda K_n)$. 

For the upper bound, let us consider $\lambda K_n=(1+t)K_n=K_n \cup tK_n$. We colour $K_n$ as in Theorem \ref{theo_main} with $\left\lceil \frac{3n}{2}\right\rceil$ colours. If $n$ is even, $tK_n$ accepts a colouring by perfect matchings with $t(n-1)$ classes of new colours, while if $n$ is odd, $tK_n$ accepts a colouring by maximal matchings with $tn$ classes of new colours. For the exceptional case of $n=4$, it is not hard to see that $\mathcal{G}^*$ is a subgraph of $K_{1,3}$ or the multigraph $K_1\cup K_1 \cup 2K_2$, and in any case its complement contains a triangle which belongs to three different classes and not two that are considered in the count since each remaining class is a matching, therefore $h_t(\lambda K_{4})\geq 3\lambda+4$, and the result follows.
\end{proof}

\section*{Acknowledgments}
The results presented in this short note were achieved through collaborative efforts during the workshop ``Extremal Graphs arising from Designs and Configurations'' hosted by the Banff International Research Station (BIRS) between May 14\textsuperscript{th} and May 19\textsuperscript{th}, 2023.  The authors express gratitude to the organizers for the opportunity and the welcoming atmosphere of the workshop. This work reflects the positive and productive environment fostered by the organizers and participants alike.

\bibliographystyle{amsplain}
\bibliography{harmonious}

\providecommand{\bysame}{\leavevmode\hbox to3em{\hrulefill}\thinspace}
\providecommand{\MR}{\relax\ifhmode\unskip\space\fi MR }
\providecommand{\MRhref}[2]{%
  \href{http://www.ams.org/mathscinet-getitem?mr=#1}{#2}
}
\providecommand{\href}[2]{#2}
\begin{thebibliography}{10}

\bibitem{AAEEJR}
A.~Aflaki, S.~Akbari, K.J. Edwards, D.S. Eskandani, M.~Jamaali, and
  M.~Ravanbod, \emph{On harmonious colourings of trees}, Electron. J. Combin.
  \textbf{19} (2012), 3.

\bibitem{AKK}
S.~Akbari, J.~Kim, and A.~Kostochka, \emph{Harmonious coloring of trees with
  large maximum degree}, Discret Math. \textbf{312} (2012), 1633--1637.

\bibitem{AMOR}
G.~Araujo-Pardo, Juan~Jos\'e Montellano-Ballesteros, Mika Olsen, and Christian
  Rubio-Montiel, \emph{On the harmonious chromatic number of graphs}, 2023.

\bibitem{BMNRM}
M.~Buratti, F.~Merola, A.~Naki\'c, and C.~Rubio-Montiel, \emph{Banff designs:
  difference methods for coloring incidence graphs}, submitted.

\bibitem{EMcD}
K.J. Edwards and C.J.H. McDiarmid, \emph{New upper bounds on harmonious
  colorings}, J. Graph Theory \textbf{18} (1994), 257--267.

\bibitem{FHP_l_distinguishing_number}
O.~Frank, F.~Harary, and M.~Plantholt, \emph{The line-distinguishing chromatic
  number of a graph}, Ars Combin. \textbf{14} (1982), 241--252.

\bibitem{HK_Harm_coloring_of_graphs}
J.E. Hopcroft and M.S. Krishnamoorthy, \emph{On the harmonious coloring of
  graphs}, SIAM J. Algebraic Discrete Methods \textbf{4} (1983), 306--311.

\bibitem{Lu}
Z.~Lu, \emph{On an upper bound for the harmonious chromatic number of a graph},
  J. Graph Theory \textbf{15} (1991), 345--347.

\bibitem{MP_harmonious_cn}
Z.~Miller and D.~Pritikin, \emph{The harmonious coloring number of a graph},
  Discrete Math. \textbf{93} (1991), 211--228.

\bibitem{Mitchem}
J.~Mitchem, \emph{On the harmonious chromatic number of a graph}, Discrete
  Math. \textbf{74} (1989), 151--157.

\bibitem{Stong}
R.A. Stong, \emph{On $1$-factorizability of cayley graphs}, J. Comb. Theory
  Ser. B \textbf{39} (1985), 298--307.

\bibitem{VVK}
M.~Venkatachalam, J.V. Vivin, and K.~Kaliraj, \emph{Harmonious coloring of
  double star graph families}, Tamkang J. Math. \textbf{43} (2012), 153--158.

\end{thebibliography}

\end{document}